\documentclass[24pt,reqno]{article}

\usepackage{amsmath,amssymb,amsfonts,enumerate,amsthm,graphicx,color}
\usepackage{amsmath}
\usepackage{amssymb}
\usepackage{color}
\usepackage{epsf}
\usepackage{fleqn}

\usepackage{subfig}

\usepackage{amssymb}
\usepackage{amsmath}
\usepackage{amssymb}
\usepackage{amsmath}
\usepackage[latin1]{inputenc}
\usepackage[english]{babel}

\usepackage[english]{babel}
\usepackage{pgf,pgfarrows,pgfnodes,pgfautomata,pgfheaps}
\usepackage{amsmath,amssymb}
\usepackage[latin1]{inputenc}
\usepackage{graphicx,amsmath,amsfonts,amssymb,amsthm,mathrsfs}

\begin{document}
\makeatletter

\begin{center}
\epsfxsize=10in
\end{center}

\def\endofproofmark{$\Box$}

\def\endofproofmark{$\Box$}
\begin{center}
\vskip 1cm {\LARGE\bf Some Monotonicity Properties of Convex}\\ \vspace{0.5cm} {\LARGE\bf Functions with Applications}
\vskip 1cm \large
Jamal Rooin and Hossein Dehghan
\\
\vskip .5cm
Department of Mathematics\\
Institute for Advanced
Studies in Basic Sciences\\
Zanjan 45137-66731, Iran\\
\textcolor[rgb]{1.00,0.00,1.00}{{\tt rooin@iasbs.ac.ir}}\\ \textcolor[rgb]{1.00,0.00,1.00}{{\tt hossein.dehgan@gmail.com}}\\
\end{center}

\date{}

%%%%%%%%%%%%%%%%%%%%%%%%%MACROS%%%%%%%%%%%%%%%%%%%%%%%%%%%%%%%%%%%%%%%%%%%
\newtheorem{theo}{Theorem}[section]
\newtheorem{prop}[theo]{Proposition}
\newtheorem{rem}[theo]{Remark}
\newtheorem{lemma}[theo]{Lemma}
\newtheorem{cor}[theo]{Corollary}
\newtheorem{problem}{Problem}
\numberwithin{equation}{section}

\def\frameqed{\framebox(5.2,6.2){}}
\def\deshqed{\dashbox{2.71}(3.5,9.00){}}
\def\ruleqed{\rule{5.25\unitlength}{9.75\unitlength}}
\def\myqed{\rule{8.00\unitlength}{12.00\unitlength}}
\def\qed{\hbox{\hskip 6pt\vrule width 7pt height11pt depth1pt\hskip 3pt}
\bigskip}
\newcommand{\COM}[2]{{#1\choose#2}}
%%%%%%%%%%%%%%%%%%%END%OF%MACROS%%%%%%%%%%%%%%%%%%%%%%%%%%%%%%%%%%%%%%%%%%%%%%
\thispagestyle{empty} \null \addtolength{\textheight}{1cm}

\begin{abstract}
We mainly establish a monotonicity property between some special Riemann sums of a convex function $f$ on $[a,b]$, which in particular yields that $\frac{b-a}{n+1}\sum_{i=0}^n
f\left(a+i\frac{b-a}{n}\right)$ is decreasing while $\frac{b-a}{n-1}\sum_{i=1}^{n-1}
f\left(a+i\frac{b-a}{n}\right)$ is an increasing sequence. These give us a new refinement of the Hermitt-Hadamard inequality.
 Moreover, we give a refinement of the classical
Alzer's inequality together with a suitable converse to it. Applications regarding to some important convex functions are also included.
\end{abstract}

\bigskip
\hrule
\bigskip

\noindent 2010 {\it Mathematics Subject Classification}: 26D15, 26A51, 26A06.

\noindent \emph{Keywords and phrases: Convexity, Mean, Hermitt-Hadamard inequality, Alzer inequality, Bennett inequality.}

\bigskip
\hrule
\bigskip
\section{Introduction}

Throughout this paper, we suppose that $a<b$ are two real numbers and
%%$f:[a,b]\rightarrow\mathbb{R}$
$f$ is a real-valued
function on the closed interval $[a,b]$. We put
%\begin{eqnarray*}\label{3}
%M_n=\sum_{i=0}^n f\left(x_i^{(n)}\right),\hspace{1cm}
%M'_n=M_n-f(a)-f(b)=\sum_{i=1}^{n-1} f\left(x_i^{(n)}\right),
%\end{eqnarray*}
\begin{eqnarray*}
A_n:= \frac{b-a}{n+1}\sum_{i=0}^n f\left(x_i^{(n)}\right), \hspace{2cm}
B_n:=\frac{b-a}{n-1}\sum_{i=1}^{n-1} f\left(x_i^{(n)}\right),%\hspace{1cm}(n=2,3,\ldots)
\end{eqnarray*}
and
\begin{eqnarray*}\label{4} S_n=\frac{b-a}{n}\sum_{i=1}^n
f\left(x_i^{(n)}\right),\hspace{1cm}
T_n=\frac{b-a}{n}\sum_{i=0}^{n-1} f\left(x_i^{(n)}\right),
\end{eqnarray*}
where
\begin{eqnarray*}
x_i^{(n)}=a+i\frac{b-a}{n}\hspace{1.5cm}(i=0,1,\ldots,n;~n=1,2,\ldots).
\end{eqnarray*}
($B_n$ is defined for $n\geq 2$.) %When emphasizing, we write $A_n(f)$ instead of $A_n$, and so on.
\\
It is known (see, e.g., \cite[p. 565]{ABP}) that if $f$ is increasing convex or increasing concave on $[0,1]$, then the
%sequence
sequence $S_n$ is decreasing while $T_n$ is increasing, i.e.
\begin{eqnarray}\label{Sn Tn}
S_{n+1}\leq S_n,\hspace{1.5cm} T_n\leq T_{n+1}.%\hspace{1.5cm}(i=0,1,\ldots,n;~n=1,2,\ldots).
\end{eqnarray}
The inequalities in (\ref{Sn Tn}) are strict if $f$ is strictly increasing and convex or strictly increasing and concave.
\par
%In 1964. H. Minc and L. Sathre in [24] proved that, for n ¸ N the inequality
In 1964, H. Minc and L. Sathre \cite{9Min} proved that
%the following inequality.  If we let $r\rightarrow 0+$, we get the following inequality, due to H. Minc and L. Sathre \cite{9Min}:
\begin{eqnarray}\label{2}
\frac{n}{n+1}\leq\frac{\sqrt[n]{n!}}{\sqrt[n+1]{(n+1)!}}\hspace{1.5cm}(n=1,2,\ldots).
\end{eqnarray}
In 1988, J.S. Martins \cite{8Mar} established that
%, gave another lower bound for the ratio $\frac{\sqrt[n]{n!}}{\sqrt[n+1]{(n+1)!}}$
 for each $r>0$, %and positive integer $n$, %we have
\begin{eqnarray}\label{martin}
\left((n+1)\sum_{i=1}^ni^r\left/ n\sum_{i=1}^{n+1}i^r\right.\right)^{1/r} \leq \frac{\sqrt[n]{n!}}{\sqrt[n+1]{(n+1)!}}\hspace{1.5cm}(n=1,2,\ldots).
\end{eqnarray}
In 1992, G. Bennett % Lemma 8 in Bennett
\cite{Bennett} proved the following inequality
\begin{eqnarray}\label{bennett}
\left((n+1)\sum_{i=1}^ni^r\left/ n\sum_{i=1}^{n+1}i^r\right.\right)^{1/r}  \leq \frac{n+1}{n+2} \hspace{1.5cm} ( r>1;\ n=1,2,\ldots),
\end{eqnarray}
which is reversed if $r<1$.\\
In 1993, H. Alzer \cite{2Al} came into comparing the left-hand sides of (\ref{2}) and (\ref{martin}) and proved that  for each $r>0$,
\begin{eqnarray}\label{1}
\frac{n}{n+1}\le\left((n+1)\sum_{i=1}^ni^r\left/ n\sum_{i=1}^{n+1}i^r\right.
\right)^{1/r}\hspace{1.5cm}(n=1,2,\ldots).
\end{eqnarray}
%In other words, the Riemann sums $\frac{1}{n}\sum_{i=1}^n\left(\frac{i}{n}\right)^r~(n=1,2,\ldots)$ of the function $x^r$ form a decreasing sequence.
The proof of Alzer is technical, but quite complicated. So, in several
articles Alzer's proof has been simplified, and also in many
others, this
inequality has been extended; see e.g. \cite{Ch-Drag,5CQ,7EP,10Qi,11San}, and see also \cite{ABP} for some historical notes. \\
 Obviously, the Alzer inequality (\ref{1}) and Martins inequality (\ref{martin}) simultaneously give us a refinement of Minc-Sathre inequality (\ref{2}).\\
Note that if $r\rightarrow 0+$ in (\ref{1}), we get (\ref{2}) without appealing to (\ref{martin}). \\
%So, the Alzer result is stronger than the Minc-Sathre one.
Clearly, for $r>1$ the Alzer inequality (\ref{1}) gives us a reverse of Bennet inequality (\ref{bennett}), while as considering $n/(n+1)< (n+1)/(n+2)$, the Bennet inequality for $0<r<1$ is a refinement of Alzer inequality.\\
In 1994, H. Alzer \cite{Alzer1994} showed that if $r<0$, the Martins inequality (\ref{martin}) is reversed. This result is reobtained by C.P. Chen et al. \cite{Chen-Qi-Dragomir} in 2005, too.\\
Recently, J. Rooin et al. \cite{RAM}, using some technics of convexity, generalized the Alzer and Bennett inequalities to operators when $-1\leq r\leq 2$.
\par
%J. van de Lint proved that if $f$ is increasing convex or increasing concave on $[0,1]$,
Let $f$ be convex on $[a,b]$.
%In this paper, First, without any monotonicity assumptions on $f$, we prove that the sequence $A_n$ is decreasing while $B_n$ is increasing.\\
The main purpose of this paper is to prove the inequality  (\ref{new22}) regarding some Riemann sums of $f$. This inequality yields that the sequence $A_n$ is decreasing while $B_n$ is increasing, without any monotonicity assumptions on $f$. As a consequence, we give an extension and a refinement to the well-known Hermitt-Hadamard inequality \cite{MPF}:
\begin{eqnarray}\label{10}
f\left(\frac{a+b}{2}\right)\leq \frac{1}{b-a} \int _a^b f(t)dt\leq
\frac{f(a)+f(b)}{2}.
\end{eqnarray}
For more details see  \cite{6DP}. If in addition $f$ is increasing, we get some refinements and converses to (\ref{Sn Tn}).\\
%Then, applying these results for some important convex functions, we give several applications.
Applying these results to the power function $x^r$, we get the Bennett inequality (\ref{bennett}) and refinements and converses of the classical Alzer inequality (\ref{1}) in the case of $-\infty< r < +\infty$. These extend the numerical results of \cite{RAM}. Also, we obtain new inequalities concerning $p$-logarithmic means. Finally, we give applications regarding to some other important convex functions, which in particular, yield us new rational approximations of trigonometric functions.
%In particular, we obtain some refinements and converses to Alzer inequality all inequalities involved in the above table except than Martins inequality for the cases $0<r<1$ and $r<0$. In this article, first using some trivial facts about convex functions, we obtain some valuable results concerning two important special kinds of modified Riemann sums of convex functions. Then, applying these results for some important convex functions, we give several applications. In particular, we obtain all inequalities involved in the above table except than Martins inequality for the cases $0<r<1$ and $r<0$.In particular, if $f$ is increasing, these refine and give converses to (\ref{Sn Tn}). yield  refinements and  In this article, first using some trivial facts about convex functions, we obtain some valuable results concerning two important special kinds of modified Riemann sums of convex functions. Then, applying these results for some important convex functions, we give several applications. In particular, we obtain all inequalities involved in the above table except than Martins inequality for the cases $0<r<1$ and $r<0$.

\section{ Main results}

%All we need is the following theorem which gives us two recursive inequalities concerning the symmetric sums $M_n$'s and $M'_n$'s which we call them in this paper; modified Riemann sums.
In this section, we prove some monotonicity properties of convex functions. The following theorem is the main source of all results in this paper.

\begin{theo}\label{newmain}
Let $a=x_0< x_1< \cdots < x_n=b$ and $a=y_0< y_1< \cdots < y_{n+1}=b$ be two partitions of $[a,b]$ such that $x_{i-1}\leq y_i\leq x_{i}$ ($i=1,2,\ldots, n$). If $f$ is convex on $[a,b]$, then
\begin{eqnarray}\label{new2}
\sum_{i=1}^{n} (x_{i}-x_{i-1}) f(y_i) \leq \sum_{i=0}^{n} (x_{i+1}-x_{i-1}+ y_{i}-y_{i+1}) f(x_i)
\end{eqnarray}
and
\begin{eqnarray}\label{new1}
\sum_{i=0}^{n} (y_{i+1}-y_i) f(x_i) \leq \sum_{i=0}^{n+1} (y_{i+1}-y_{i-1}+ x_{i-1}-x_i) f(y_i),
\end{eqnarray}
where $x_{-1}=y_{-1}=a$ and $x_{n+1}=y_{n+2}=b$.
If $f$ is strictly convex, then inequality (\ref{new2}) (respectively (\ref{new1})) is  strict whenever
$x_{i-1}< y_i< x_{i}$ for some $i\in \{1,2,\ldots, n\}$
(respectively $y_{i}< x_i< y_{i+1}$ for some $i\in \{1,2,\ldots, n-1\}$).
%If $f$ is strictly convex, then inequality (\ref{new1}) (respectively (\ref{new2})) is  strict whenever $y_{i}< x_i< y_{i+1}$ (respectively $x_{i-1}< y_i< x_{i}$) for some $i\in \{1,2,\ldots, n-1\}$ (respectively $i\in \{1,2,\ldots, n\}$).
\end{theo}
\noindent{\it Proof}.
Since $x_{i-1}\leq y_i\leq x_{i}$ ($i=1,2,\ldots, n$), using
$$y_i= \frac{x_i-y_i}{x_{i}-x_{i-1}}\ x_{i-1} + \frac{y_{i}-x_{i-1}}{x_{i}-x_{i-1}}\ x_i$$
  and convexity of $f$ we have
\begin{eqnarray}\label{new3}
(x_{i}-x_{i-1}) f(y_i) \leq (x_i-y_i) f(x_{i-1}) + (y_{i}-x_{i-1}) f(x_i)\hspace{1.5cm} (i=1,2,\ldots, n).
\end{eqnarray}
Now summing up (\ref{new3}) from $1$ to $n$, we get
\begin{eqnarray}\label{new4}
\nonumber \sum_{i=1}^{n} (x_{i}-x_{i-1}) f(y_i) &\leq& \sum_{i=1}^{n} (x_i-y_i) f(x_{i-1}) + \sum_{i=1}^{n} (y_{i}-x_{i-1}) f(x_i)\\
\nonumber &=&\sum_{i=0}^{n-1} (x_{i+1}-y_{i+1}) f(x_{i}) + \sum_{i=1}^{n} (y_{i}-x_{i-1}) f(x_i)\\
\nonumber &=&\sum_{i=0}^{n} (x_{i+1}-y_{i+1}) f(x_{i}) + \sum_{i=0}^{n} (y_{i}-x_{i-1}) f(x_i)\\
\nonumber &=&\sum_{i=0}^{n} (x_{i+1}-x_{i-1}+ y_{i}-y_{i+1}) f(x_i).
\end{eqnarray}
%and (\ref{new1}) holds.
The inequality (\ref{new1}) follows in a similar manner by considering  $y_{i}\leq x_{i}\leq y_{i+1}$ ($i=0,1,\ldots, n$). The rest is clear. $\Box$
\begin{rem}
With the assumptions of Theorem \ref{newmain} we may write the inequalities (\ref{new2}) and (\ref{new1}) in following single form
\begin{align}\label{new22}
\nonumber \sum_{i=0}^{n} &(y_{i+1}-y_i) f(x_i) +  \sum_{i=0}^{n-1} (x_{i+1}-x_{i}) f(y_{i+1}) \\
&\leq \min\left\{ \sum_{i=0}^{n} (y_{i+1}-y_i) (f(y_i)+f(y_{i+1})) , \sum_{i=0}^{n-1} (x_{i+1}-x_{i}) (f(x_i)+ f(x_{i+1}))\right\},
\end{align}
which is a monotonicity property between some special Riemann sums.
\end{rem}
%%%%%%%%%%%%%%%%%%%%%%%%%%%%%%%%%%%%%%%%%%%%%%%%%%%%%%%%%%%%
\begin{cor}\label{main}
  With the above assumptions, if $f$ is
convex on $[a,b]$, then we have
\begin{eqnarray}\label{5}
A_{n+1}\leq A_n\hspace{1cm} (n=1,2,\ldots)\hspace{1cm} \mbox{and} \hspace{1cm}
B_n\leq B_{n+1}\hspace{1cm} (n=2,3,\ldots).
\end{eqnarray}
Both inequalities are
strict if $f$ is strictly convex.
%If $f$ is concave, all inequalities reverse.
\end{cor}
\noindent{\it Proof}.
Take $x_i=x_i^{(n)}$ ($i=0,1,\ldots, n$) and $y_i=x_i^{(n+1)}$ ($i=0,1,\ldots, n+1$) in Theorem \ref{newmain}.
$\Box$
%%%%%%%%%%%%%%%%%%%%%%%%%%%%%%%%%%%%%%%%%%%%%%%%%%%%%%%%%%%%%%%%%%%%%%%%%%%%%%%%%%%%%%%%%%%%%%%%%%%%%%%%%%%%%%%%%%%%%%%%%%%%%
\begin{cor} If $f$ is convex on $[a,b]$, then
\begin{eqnarray}\label{9}
\hspace{-0.75cm}\frac{1}{m-1}\sum_{i=1}^{m-1} f\left(x_i^{(m)}\right)\leq
\frac{1}{b-a}\int _a^b f(t)dt\leq \frac{1}{n+1}\sum_{i=0}^
nf\left(x_i^{(n)}\right)\hspace{0.5cm} (m=2,3,\ldots;~n=1,2,\ldots),
\end{eqnarray}
which is a refinement and extension of Hermitt-Hadamard inequality (\ref{10}).\\
Both inequalities in (\ref{9}) are strict, if $f$ is strictly convex.
 \end{cor}
\noindent{\it Proof}.
%If $f$ be a convex function on $[a,b]$.
Clearly $f$ is Riemann integrable on $[a,b]$ and
\begin{eqnarray}\label{Limit}
\lim_{n\to\infty} A_n =\lim_{n\to\infty} B_n= \int _a^b f(t)dt.
\end{eqnarray}
Now, (\ref{9}) follows from Corollary \ref{main}.
$\Box$
%{\it With the above assumptions, if $f$ is
%(strictly) convex, then the arithmetic means
%$$A_n:=\frac{1}{n+1}\sum_{i=0}^n
%f\left(x_i^{(n)}\right)(strictly)\hspace{1cm}(n=1,2,\cdots),$$ and
%$$A'_n:=\frac{1}{n-1}\sum_{i=1}^{n-1}
%f\left(x_i^{(n)}\right)$$
%form (strictly) increasing and (strictly) decreasing sequences respectively.
%%%%%%%%%%%%%%%%%%%%%%%%%%%%%%%%%%%%%%%%%%%%%%%%%%%%%%%%%%%%%%%%%%%%%%%%%%%%%%%%%%%%%%%%%%%%%%%%%%%%%%%%%%%%%%%%%%%%%%
\begin{cor}  If $f>0$ is logarithmically convex on
$[a,b]$, then
\begin{eqnarray}\label{11}
\frac{\sqrt[m-1]{\prod_{i=1}^{m-1}f\left(\frac{(m-i)a+ib}
{m}\right)}}{\sqrt[m]{\prod_{i=1}^mf\left(\frac{(m+1-i)a+ib}
{m+1}\right)}}\leq1\leq
\frac{\sqrt[n+1]{\prod_{i=0}^{n}f\left(\frac{(n-i)a+ib}
{n}\right)}}{\sqrt[n+2]{\prod_{i=0}^{n+1}f\left(\frac{(n+1-i)a+ib}
{n+1}\right)}}
\end{eqnarray}
and
\begin{eqnarray}\label{12}
\nonumber\sqrt[m-1]{\prod_{i=1}^{m-1}f\left(\frac{(m-i)a+ib}
{m}\right)}&\leq& \exp\left(\frac{1}{b-a}\int _a^b \ln f(t)dt\right)\\
&\leq&
\sqrt[n+1]{\prod_{i=0}^{n}f\left(\frac{(n-i)a+ib} {n}\right)},
\end{eqnarray}
where
$m=2,3,\ldots$ and $n=1,2,\ldots$. \\
All inequalities are strict if $f$ is strictly logarithmically convex.
\end{cor}
\noindent{\it Proof}. Take $\ln f$ instead of $f$ in (\ref{5}) and (\ref{9}). \ $\Box$
%%%%%%%%%%%%%%%%%%%%%%%%%%%%%%%%%%%%%%%%%%%%%%%%%%%%%%%%%%%%%%%%%%%%%%%%%%%%%%%%%%%%%%%%%%%%%%%%%%%%%%%%%%%%%%%%%%%%%%%%%%%
\begin{cor}  With the above assumptions, if $f$ is
convex on $[a,b]$, then we have
\begin{eqnarray}\label{13'}
 \frac{1}{n(n+2)}\left[S_{n+1}-(b-a)f(a))\right]\le S_n-S_{n+1}\le \frac{1}{n^2}\left[(b-a)f(b)-S_{n+1}\right] \end{eqnarray}
 and
 \begin{eqnarray}\label{14'}
  \frac{1}{n^2}\left[T_{n+1}-(b-a)f(a)\right] \le  T_{n+1}-T_n\le \frac{1}{n(n+2)}\left[(b-a)f(b)-T_{n+1}\right] .
 \end{eqnarray}
%\begin{eqnarray}\label{13} S_{n+1}+\frac{1}{n(n+2)}\left[S_{n+1}-(b-a)f(a))\right]\le S_n\le S_{n+1}+\frac{1}{n^2}\left[(b-a)f(b)-S_{n+1}\right] \end{eqnarray} and \begin{eqnarray}\label{14} T_{n+1}+\frac{1}{n(n+2)}\left[T_{n+1}-(b-a)f(b))\right]\le T_n\le T_{n+1}+\frac{1}{n^2}\left[(b-a)f(a)-T_{n+1}\right]. \end{eqnarray}
Moreover, except than the case $n=1$ in which equality always
holds in the right of (\ref{13'}) and left hand of (\ref{14'}), all inequalities are
strict if $f$ is strictly convex.\\
If $f$ is concave, all inequalities reverse.
\end{cor}
\noindent{\it Proof}.
The left inequality of (\ref{13'}) and right inequality of (\ref{14'}) follow from the left hand of (\ref{5}),
% in Theorem \ref{main},  by decreasing
 %(\ref{6})
 by considering $$A_{n}=\frac{n}{n+1}S_{n}+ \frac{b-a}{n+1}f(a)\hspace{.5cm} \mbox{and} \hspace{.5cm}A_{n}=\frac{n}{n+1}T_{n}+ \frac{b-a}{n+1}f(b)\hspace{1.5cm} (n=1,2,\ldots).$$ Obviously, equality holds in right hand of (\ref{13'}) and left hand of (\ref{14'}) if $n=1$.
%Let $n\geq 2$.
Now, if $n\geq 2$, the right inequality of (\ref{13'}) and the left inequality of (\ref{14'}) follow from the right hand of (\ref{5}),
% in Theorem \ref{main},
by considering
% the fact that $B_n$ is increasing and
% (\ref{6}) by considering
$$B_{n}=\frac{n}{n-1}S_{n}- \frac{b-a}{n-1}f(b)\hspace{.5cm} \mbox{and} \hspace{.5cm}B_{n}=\frac{n}{n-1}T_{n}- \frac{b-a}{n-1}f(a)\hspace{1.5cm} (n=2,3,\ldots).$$
%Finally, inequalities (\ref{15}) and (\ref{16}) are trivially obtained by comparing the left and right hand sides of (\ref{13'}) and (\ref{14'}) with each other, respectively.\\
If $f$ is strictly convex, the strictness of all inequalities follow from strictness of inequalities in (\ref{5}).\\
%The last assertion is trivial.
$\Box$
%%%%%%%%%%%%%%%%%%%%%%%%%%%%%%%%%%%%%%%%%%%%%%%%%%%%%%%%%%%%%%%%%%%%%%%%%%%%%%%%%%%%%%%%%%%%%%%%%%%%%%%%%%%%%%%%%%%%%%%%%
\begin{rem}
%(See \cite{ABP})
If $f$ is increasing and convex (concave) on $[a,b]$,
%since
 the inequalities in (\ref{13'}) and (\ref{14'}) (the reversed forms of the inequalities in (\ref{13'}) and (\ref{14'})) give us a refinement and converse to the inequalities in (\ref{Sn Tn}).
%or increasing concave  function on $[a,b]$,since  \begin{eqnarray}\label{18new} (b-a)f(a)\leq S_{n+1}\leq(b-a)f(b) \end{eqnarray} and \begin{eqnarray}\label{19new} (b-a)f(a)\leq T_{n+1}\leq(b-a)f(b), \end{eqnarray}
%The next corollary shows that the inequalities  by (\ref{13}) and (\ref{14}) we have
 %are actually give us refinements and converses to the known inequalities (\ref{Sn Tn}).\\ \begin{eqnarray}\label{13'} 0\leq \frac{1}{n(n+2)}\left[S_{n+1}-(b-a)f(a))\right]\le S_n-S_{n+1}\le \frac{1}{n^2}\left[(b-a)f(b)-S_{n+1}\right] \end{eqnarray} and \begin{eqnarray}\label{14'} 0\leq \frac{1}{n^2}\left[T_{n+1}-(b-a)f(a)\right] \le  T_{n+1}-T_n\le \frac{1}{n(n+2)}\left[(b-a)f(b)-T_{n+1}\right] . \end{eqnarray} Except than the case $n=1$ in which equality always holds in the last and second inequalities in (\ref{13}) and (\ref{14}) respectively, all inequalities are strict if $f$ is strictly increasing.\\ These give us refinements and converses to (\ref{Sn Tn}). If $f$ is increasing concave, similar results hold. We omit the details.
 %all inequalities reverse.\\ Similarly, if $f$ is increasing concave
%or increasing concave function on $[a,b]$, the revered forms of inequalities (\ref{13}), (\ref{14}), (\ref{18new}) and (\ref{19new}) give us refinements and converses to (\ref{Sn Tn}). \begin{eqnarray}\label{17} S_{n+1}\leq S_n\hspace{1cm}\mbox{and}\hspace{1cm}T_n\leq T_{n+1}\hspace{1cm}(n=1,2,\ldots). \end{eqnarray} Both inequalities in (\ref{17}) are strict, if $f$ is strictly convex or concave.
\end{rem}
%\noindent{\it Proof}. Since $f$ is increasing, for each $k$, we have \begin{eqnarray}\label{18} (b-a)f(a)\leq S_k\leq(b-a)f(b) \end{eqnarray} and \begin{eqnarray}\label{19} (b-a)f(a)\leq T_k\leq(b-a)f(b). \end{eqnarray} So, if $f$ is convex, using (\ref{18}), (\ref{19}), the left hand side of (\ref{13}) and the right hand side of (\ref{14}), we get (\ref{17}).\\ Now, if $f$ is concave, then $-f$ is convex, and (\ref{17}) follows from the right hand side of (\ref{13}) and the left hand side of (\ref{14}) applying to $-f$, by using (\ref{18}) and (\ref{19}) and taking into account that $$ S_k(-f)=-S_k(f)\hspace{1cm}\mbox{and}\hspace{1cm}T_k(-f)=-T_k(f)\hspace{1cm}(k=1,2,\ldots). $$ In the case of strict convexity or concavity, the strictness of the inequalities in (\ref{17}) follow from the strictness of the inequalities in (\ref{13}) and (\ref{14}). (The cases $S_2<S_1$ and $T_1<T_2$ must be treated separately.) $\Box$
\section{Applications}
In this section, using the results of the preceding one, we give several nice applications
regarding some important convex functions.
\subsection{Applications to normed spaces\\}
%\begin{rem}
Let $X$ be a real normed linear space,
$x,y\in X$ and $p\geq 1$. It is clear that
$$
\varphi(t)= \|(1-t)x+ty\|^p\hspace{1.5cm}(t\in \mathbb{R})$$ is a
convex function on the real line. If $X$ is strictly convex and
$x,y$ are linearly independent, then using $\|u+v\|<\|u\|+\|v\|$ for any linearly independent vectors $u$ and
$v$, we see that $t\rightarrow \|(1-t)x+ty\|$ is strictly convex on
$\mathbb{R}$. Now, since the function $t\rightarrow t^p$ is convex and strictly increasing on $[0,\infty)$, we conclude that $\varphi$ is strictly convex on
$\mathbb{R}$.
%\end{rem}
\begin{theo} Let $x,y$ be
two vectors in a real normed linear space $X$, not both of them zero, and $p\geq 1$. Then
\begin{eqnarray*}\label{36}
\hspace{-0.5cm}\left(\frac{n\sum_{i=1}^{n-1}\|(n-i)x+iy\|^p}{(n-1)\sum_{i=1}^n\|(n+1-i)x+iy\|^p}\right)^{1/p}
\leq\frac{n}{n+1}\leq\left(\frac{(n+2)\sum_{i=0}^n\|(n-i)x+iy\|^p}{(n+1)\sum_{i=0}^{n+1}\|
(n+1-i)x+iy\|^p}\right)^{1/p}
\end{eqnarray*}
and
\begin{eqnarray}\label{37}
\frac{\sum_{i=1}^{n-1}\|(n-i)x+iy\|^p}{n^p(n-1)} \leq
\int_0^1\|(1-t)x+ty\|^pdt\leq\frac{\sum_{i=0}^n\|(n-i)x+iy\|^p}{n^p(n+1)},
\end{eqnarray}
where in the left hands $n\geq2$ and in the right hands $n\geq1$.\\
Note that (\ref{37}) is a generalization and refinement of the well-known chain inequalities \cite{MPF}
\begin{eqnarray*}\label{38}
\left\|\frac{x+y}{2}\right\|^p\leq
\int_0^1\|(1-t)x+ty\|^pdt\leq\frac{\|x\|^p+\|y\|^p}{2}.
\end{eqnarray*}
If $X$ is strictly convex, then all inequalities are strict if
$x$ and $y$ are linearly independent.
\end{theo}
\noindent{\it Proof}. Apply (\ref{5}) and (\ref{9}) to the convex function $\varphi$ on $[0,1]$. $\Box$
%$$
%M_k=\sum_{i=0}^k\|(1-\frac{i}{k})x+\frac{i}{k}y\|^p=\frac{1}{k^p}\sum_{i=0}^k\|(k-i)x+iy\|^p\hspace{1.5cm}(k=1,2,\cdots).
%$$
%Now, (\ref{21}) follows from (\ref{5}).\\\\
\subsection{ Applications to power and logarithmic functions\\}
We recall that the $p$-logarithmic, identric and logarithmic means
of $a,b> 0$ are defined respectively by
\begin{eqnarray*}\label{20} L_p(a,b)= \left\{
\begin{array}{cl}
a &\mbox{if}\hspace{2mm}a=b\\
\left[\frac{b^{p+1}-a^{p+1}}{(p+1)(b-a)}\right]^{1/p}&\mbox{if}\hspace{2mm}a\not=b
\end{array}
\right.,\hspace{1.5cm} p\in \mathbb{R}\setminus \{0,-1\},
\end{eqnarray*}
\begin{eqnarray*}
I(a,b)= \left\{
\begin{array}{cl}
a &{\rm if\hspace{2mm}}a=b\\
\frac{1}{e}\left(\frac{b^b}{a^a}\right)^\frac{1}{b-a} &{\rm
if}\hspace{2mm}a\not=b
\end{array}
\right.
\end{eqnarray*}
and
\begin{eqnarray*}\label{22}
L(a,b)= \left\{
\begin{array}{cl}
a &\mbox{if}\hspace{2mm}a=b\\
\frac{b-a}{\ln b-\ln a}&\mbox{if}\hspace{2mm}a\not=b
\end{array}
\right. .
\end{eqnarray*}
Note that
\begin{eqnarray}\label{23} \lim_{p\rightarrow
0}L_p(a,b)=I(a,b)\hspace{1.5cm}\mbox{and}\hspace{1.5cm}
\lim_{p\rightarrow -1}L_p(a,b)=L(a,b).\nonumber
\end{eqnarray}
So, we can take $L_0=I$ and $L_{-1}=L$. Note that $L_p(a,b)$ is also defined if $0\not=p>-1$ and $a,b\geq 0$.
\begin{theo} Let $0\leq a<b$. If $r>1$, then
\begin{eqnarray}\label{24}
\hspace{-0.5cm}\left(\frac{n\sum_{i=1}^{n-1}[(n-i)a+ib]^r}{(n-1)\sum_{i=1}^n[(n+1-i)a+ib]^r}\right)^{1/r}
<\frac{n}{n+1}<\left(\frac{(n+2)\sum_{i=0}^n[(n-i)a+ib]^r}{(n+1)\sum_{i=0}^{n+1}[(n+1-i)a+ib]^r}\right)^{1/r}
\end{eqnarray}
and
\begin{eqnarray}\label{25}
\left(\frac{\sum_{i=1}^{n-1}[(n-i)a+ib]^r}{n^r(n-1)}\right)^{1/r}
<L_r(a,b)<\left(\frac{\sum_{i=0}^n[(n-i)a+ib]^r}{n^r(n+1)}\right)^{1/r},
\end{eqnarray}
where in the left hand inequalities, we have $n\geq2$, and in the
right hand ones, $n\geq1$.
\\ If $r<0$ with $a>0$ or $0<r<1$, all
inequalities in (\ref{24}) and (\ref{25}) reverse.\\
%In particular, \begin{eqnarray*}\label{26} (n+1)\left/\sum_{i=0}^n\frac{n}{(n-i)a+ib}\right.< L(a,b)<(n-1)\left/\sum_{i=1}^{n-1}\frac{n}{(n-i)a+ib}\right., \end{eqnarray*} where in the left hand, we have $n\geq1$, and in the right hand, $n\geq2$.
\end{theo}
%\begin{eqnarray}
%\frac{2}{\frac{1}{a}+\frac{1}{b}}\leq L(a,b)\leq\frac{a+b}{2}.
%\end{eqnarray}\\
\noindent{\it Proof}. For $r>1$ and $r<0$ the function $f(x)=x^r$ is strictly convex on
$[0,\infty)$ and $(0,\infty)$ respectively. So if we apply
(\ref{5}) and (\ref{9}) for $f$ on $[a,b]$, we achieve the results.\\
If $0<r<1$, the function $f$ is strictly concave on $[0,\infty)$,
and so
both inequalities in (\ref{24}) and (\ref{25}) reverse.  $\Box$
\begin{theo} If $0<a<b$, then
\begin{eqnarray}\label{26}
\frac{\sqrt[n+1]{\prod_{i=0}^n[(n-i)a+ib]}}
{\sqrt[n+2]{\prod_{i=0}^{n+1}[(n+1-i)a+ib]}}<\frac{n}{n+1}
<\frac{\sqrt[n-1]{\prod_{i=1}^{n-1}[(n-i)a+ib]}}
{\sqrt[n]{\prod_{i=1}^{n}[(n+1-i)a+ib]}}
\end{eqnarray}
and
\begin{eqnarray}\label{27}
\frac{\sqrt[n+1]{\prod_{i=0}^n[(n-i)a+ib]}} {n}<
I(a,b)<\frac{\sqrt[n-1]{\prod_{i=1}^{n-1}[(n-i)a+ib]}}{n},
\end{eqnarray}
where in the left hand inequalities, we have $n\geq1$, and in the
right hand one $n\geq2$.
\end{theo}
\noindent{\it Proof}.
Applying (\ref{11}) and (\ref{12}) for the strictly logarithmically convex function $f(x)=1/x$ on
$[a,b]$, we get (\ref{26}) and (\ref{27}). $\Box$

%Applying (\ref{6}) and (\ref{9}) for the strictly convex function $f(x)=-\ln x$ on
%$[a,b]$, we get (\ref{26}) and (\ref{27}). $\Box$
%The inequality (\ref{20}) follows from the left hand side of (\ref{19}), first
%by letting $a\rightarrow
%0+$, and then by replacing $n$ by $n+1$.
%%%%%%%%%%%%%%%%%%%%%%%%%%%%%%%%%%%%%%%%%%%%%%%%%%%%%%%%%%%%%%%%%%%%%%%%%%%%%%%%%%%%%%%%%%%%%%%%%%%%%%%%%%%%%%%%%%%%%%%%%%%%%%%%%
\begin{rem} (i) If we set $a=0$ and $b=1$ in (\ref{23}), we get
for $r>1$,
\begin{eqnarray}\label{28}
\left(\frac{n\sum_{i=1}^{n-1}i^r}{(n-1)\sum_{i=1}^n
i^r}\right)^{1/r}
<\frac{n}{n+1}<\left(\frac{(n+2)\sum_{i=1}^ni^r}{(n+1)\sum_{i=1}^{n+1}i^r}\right)^{1/r},
\end{eqnarray}
where in the left hand inequality, we have $n\geq2$, and in the
right hand one $n\geq1$. It can be seen that (\ref{28}) in turn is
equivalent to
\begin{eqnarray}\label{29}
\hspace{-0.5cm}\frac{n}{n+1}\left(1+\frac{1}{n(n+2)}\right)^{1/r}<
\left((n+1)\sum_{i=1}^ni^r\left/n\sum_{i=1}^{n+1}i^r \right.\right)^{1/r}< \frac{n+1}{n+2} \hspace{1cm}(n\geq 1).
\end{eqnarray}
%$$<\frac{n}{n+1}\left(1+\frac {(n+1)^{r+1}-\sum_{i=1}^{n+1}i^r}
%{n^2\sum_{i=1}^{n+1}i^r} \right)^{1/r},
%$$
%where $n\geq1$.\\
Similarly, if $0<r<1$, all inequalities in (\ref{28}) and so in (\ref{29})
reverse.\\
The inequalities in (\ref{29}) and their reversed forms in the case of
$0<r<1$, give us Bennett inequality (\ref{bennett}) in the case of $r>0$ and a refinement  and converse of the classical
Alzer's inequality (\ref{1}) which are stronger than the result in \cite[Corollary 1]{Ch-Drag}.\\
%%%%%%%%%%%%%%%%%%%%%%%%%%%%%%%%%%%%%%%%%%%%%%%%%%%%%%%%%%%%%%%%%%%%%%%%%%%%%%%%%%%%%%%%%%%%%%%%%%%%%%%%%%%%%%%%%%%%%%%%%%%%%%%%%
(ii) If we take $b=1$ and let $a\to 0+$ in the right hand inequality of (\ref{27}), we get
\begin{eqnarray*}\label{week Sterling}
\sqrt[n]{n!} \geq \frac{n+1}{e}  \hspace{1.5cm}(n=1,2,\ldots).
\end{eqnarray*}
\\
%%%%%%%%%%%%%%%%%%%%%%%%%%%%%%%%%%%%%%%%%%%%%%%%%%%%%%%%%%%%%%%%%%%%%%%%%%%%%%%%%%%%%%%%%%%%%%%%%%%%%%%%%%%%%%%%%%%%%%%%%%%%%%%%%
(iii) If $r<0$, letting $a\rightarrow 0+$ and changing $n$ by $n+1$ in the reversed form of the left hand inequality of (\ref{23}), we obtain
\begin{eqnarray*}\label{30}
\frac{n+1}{n+2}
\leq\left(\frac{(n+1)\sum_{i=1}^{n}i^r}{n\sum_{i=1}^{n+1}
i^r}\right)^{1/r}\hspace{1.5cm}(r<0;~n=1, 2,\ldots),
\end{eqnarray*}
which
%by changing $n$ to $n+1$, it
is the Bennett inequality (\ref{bennett}) for $r<0$.\\
%%%%%%%%%%%%%%%%%%%%%%%%%%%%%%%%%%%%%%%%%%%%%%%%%%%%%%%%%%%%%%%%%%%%%%%%%%%%%%%%%%%%%%%%%%%%%%%%%%%%%%%%%%%%%%%%%%%%%%%%%%%%%%%%%
(iv) If $r>1$, setting $a=0$ and $b=1$ in (\ref{25}), we get
\begin{eqnarray}\label{0<r<1}
 \frac{(n+1)n^r}{r+1} < \sum_{i=1}^n i^r < \frac{n(n+1)^r}{r+1}  \hspace{1.5cm}( n=1,2,\ldots).
\end{eqnarray}
If $0<r < 1$, then inequalities in (\ref{0<r<1}) reverse.\\
Also, if in the reversed form of the left hand of (\ref{25}), we take $b=1$, let $a\to 0+$ and change $n$ by $n+1$, we obtain
\begin{eqnarray*}\label{serie 1/n}
\sum_{i=1}^n i^r \leq \frac{n(n+1)^r}{r+1}  \hspace{1.5cm}(-1 < r< 0;\ n=1,2,\ldots).
\end{eqnarray*}\\
%%%%%%%%%%%%%%%%%%%%%%%%%%%%%%%%%%%%%%%%%%%%%%%%%%%%%%%%%%%%%%%%%%%%%%%%%%%%%%%%%%%%%%%%%%%%%%%%%%%%%%%%%%%%%%%%%%%%%%%%%%%%%%%%%
%%%%%%%%%%%%%%%%%%%%%%%%%%%%%%%%%%%%%%%%%%%%%%%%%%%%%%%%%%%%%%%%%%%%%%%%%%%%%%%%%%%%%%%%%%%%%%%%%%%%%%%%%%%%%%%%%%%%%%%%%%%%%%%%%
(v) If $0<a<b$, letting $r\rightarrow 0$ in the reversed form of (\ref{24}) and (\ref{25}), we get
a weaker form of (\ref{26}) and (\ref{27}), loosing the strictness of
inequalities.\\
%%%%%%%%%%%%%%%%%%%%%%%%%%%%%%%%%%%%%%%%%%%%%%%%%%%%%%%%%%%%%%%%%%%%%%%%%%%%%%%%%%%%%%%%%%%%%%%%%%%%%%%%%%%%%%%%%%%%%%%%%%%%%%%%%
(vi) If $a\rightarrow 0+$ in (\ref{26}), changing $n$ by $n+1$, we get \cite[Lemma 2.1]{Alzer1994}
\begin{eqnarray*}\label{31}
\frac{n+1}{n+2}
\leq\frac{\sqrt[n]{n!}}{\sqrt[n+1]{(n+1)!}}\hspace{1.5cm}(n=1,2,\ldots),
\end{eqnarray*}
which is a refinement of the inequality of H. Minc and L. Sathre
(\ref{2}).
\end{rem}
%\section{Other Applications}
\begin{theo} If $0<a<b\leq \frac{1}{2}$, then
\begin{eqnarray}\label{32}
\frac{\sqrt[m-1]{\prod_{i=1}^{m-1}\frac{m-(m-i)a-i b}
{(m-i)a+ib}}}{\sqrt[m]{\prod_{i=1}^{m}\frac{(m+1)-(m+1-i)a-i b}
{(m+1-i)a+ib}}}<1<\frac{\sqrt[n+1]{\prod_{i=0}^{n}\frac{n-(n-i)a-i b}
{(n-i)a+ib}}}{\sqrt[n+2]{\prod_{i=0}^{n+1}\frac{(n+1)-(n+1-i)a-i b}
{(n+1-i)a+ib}}}
\end{eqnarray}
and
\begin{eqnarray}\label{33}
\sqrt[m-1]{\prod_{i=1}^{m-1}\frac{m-(m-i)a-i b}
{(m-i)a+ib}}<\frac{I(1-a,1-b)}{I(a,b)}<\sqrt[n+1]{\prod_{i=0}^{n}\frac{n-(n-i)a-i
b} {(n-i)a+ib}}
\end{eqnarray}
$$(m=2,3,\ldots;n=1,2,\ldots).$$
In particular,
\begin{eqnarray}\label{34}
{2m-1 \choose m-1}^\frac{1}{m-1}\leq{2m+1 \choose
m}^\frac{1}{m}\hspace{1.5cm}(m=2,3,\ldots),
\end{eqnarray}
and so ${2m+1 \choose m}^\frac{1}{m}$ is an increasing
sequence which tends to $4$. Also, we have \cite[Theorem 3.4]{13RH}
\begin{eqnarray}\label{35}
\frac{2-a-b}
{a+b}<\frac{I(1-a,1-b)}{I(a,b)}<\sqrt{\frac{(1-a)(1-b)} {a
b}}.
\end{eqnarray}
\end{theo}
\noindent{\it Proof}. The function $f(t)=\frac{1-t}{t}$ is
strictly logarithmically convex on
$(0,1/2]$. So, employing (\ref{11}) and (\ref{12}) for the function $f$ on $[a,b]$ we yield (\ref{32}) and (\ref{33}). The inequality (\ref{34}) follows from the left hand inequality in (\ref{32}) by taking $b=1/2$ and letting $a\rightarrow 0+$.
Set
\begin{eqnarray}
\nonumber u(m,a)=\sqrt[m-1]{\prod_{i=1}^{m-1}\frac{m-(m-i)a-\frac{i}{2}}
{(m-i)a+\frac{i}{2}}}\hspace{1.5cm}(m=2,3,\ldots; \ 0<a<1/2).
\end{eqnarray}
Since $B_m$ is increasing and $f$ is decreasing, $u(m,a)$ is increasing with respect to $m$ and decreasing with respect to $a$. So considering (\ref{Limit}) we get
\begin{eqnarray}
\nonumber\lim_{m\to\infty} {2m+1 \choose m}^\frac{1}{m} &=&  \lim_{m\to\infty} \lim_{a\to 0+} u(m,a) =  \sup_{m} \sup_{0<a<1/2} u(m,a)= \sup_{0<a<1/2}\sup_{m} u(m,a) \\
%\nonumber &=& \sup \{u(m,a) : m=2,3,\ldots;  \ \ 0<a<1/2\}\\
\nonumber &=&  \lim_{a\to 0+}\lim_{m\to\infty} u(m,a) =\lim_{a\to 0+} \frac{I(1-a,1/2)}{I(a,1/2)} = 4.
\end{eqnarray}
Finally, (\ref{35}) is an special case of (\ref{33}) for the choices $m=2$ and $n=1$. $\Box$

\subsection{Applications to trigonometric functions\\}

We conclude this section with the following
trigonometric estimations.
\begin{theo} If $0<x\leq\pi/2$, then
\begin{eqnarray}\label{41}
\frac{n-1}{n}\cot\frac{x}{n+1}+\frac{1}{n}\cot x< \cot\frac{x}{n}<
\frac{n+1}{n+2}\cot\frac{x}{n+1}-\frac{1}{n+2}\cot x,
\end{eqnarray}
\begin{eqnarray}\label{42}
%(\cos x \leq )
\frac{1}{n+1}(\sin x\cot\frac{x}{n}+\cos x)<\frac{\sin
x}{x} < \frac{1}{n-1}(\sin x\cot\frac{x}{n}-\cos x),%(\leq 2\cos x),
\end{eqnarray}
and in particular,
\begin{eqnarray}\label{43}
\frac{n-1}{n}\cot\frac{\pi}{2(n+1)}< \cot\frac{\pi}{2n}<
\frac{n+1}{n+2}\cot\frac{\pi}{2(n+1)}.
\end{eqnarray}\\
where in the left hands of (\ref{41}) and (\ref{43}) and in the right hand of (\ref{42}) we have $n\geq2$ and in the others $n\geq1$.
\end{theo}
\noindent{\it Proof}. The function $f(t)=\sin t$ is strictly concave on
$[0,2x]\subseteq[0,\pi]$. Now, applying (\ref{5}) and (\ref{9}) in the reversed order
to $f$, and considering
\begin{eqnarray}
\nonumber \frac{n+1}{2x}A_n=\sum _{i=0}^n\sin\frac{2ix}{n}=\frac{\sin(\frac{n+1}{n}x)\sin
x}{\sin\frac{x}{n}}= \sin^2
x\cot\frac{x}{n} + \sin x\cos x%\hspace{1.5cm}(n=1,2,\ldots)
\end{eqnarray}
and
\begin{eqnarray}\nonumber
\frac{n-1}{2x}B_n=\sum _{i=1}^{n-1}\sin\frac{2ix}{n}=\sum _{i=1}^{n}\sin\frac{2ix}{n}-\sin 2x
=\sin^2 x\cot\frac{x}{n}- \sin x\cos x
%=\frac{n-1}{2x}A_{n-1} =\sin x\cos x+\sin^2 x\cot\frac{x}{n-1} \hspace{1.5cm}(n=2,3,\ldots),
\end{eqnarray}
we obtain (\ref{41}) and (\ref{42}). The inequalities in (\ref{43}) follow from (\ref{41}) by taking $x=\pi/2$. $\Box$
\begin{rem} From (\ref{43}), we have
\begin{eqnarray}\label{44}
\frac{k-1}{k}<\frac{\tan\frac{\pi}{2(k+1)}}{\tan\frac{\pi}{2k}}<\frac{k+1}{k+2}\hspace{1.5cm}(k=2,3,\ldots),
\end{eqnarray}
which by multiplying each side of (\ref{44}) from $k=2$ to $k=n-1$, we
obtain
\begin{eqnarray}\label{45}
\frac{1}{n-1}<\tan\frac{\pi}{2n}<\frac{3}{n+1}\hspace{1.5cm}(n=3,4,\ldots).
\end{eqnarray}
Now using the representations of $\tan2x$, $\cos2x$ and
$\sin2x$ in terms of $\tan x$, and applying (\ref{45}), we get for $n=3,4,\ldots$, the following rational approximations
\begin{eqnarray}\label{46}
\frac{2(n-1)}{n(n-2)}<\tan\frac{\pi}{n}<\frac{6(n+1)}{(n+4)(n-2)},%\hspace{1.5cm}(n=3,4,\ldots),
\end{eqnarray}
\begin{eqnarray}\label{47}
\frac{(n-2)(n+4)}{n^2+2n+10}<\cos\frac{\pi}{n}<\frac{n(n-2)}{n^2-2n+2},%\hspace{1.5cm}(n=3,4,\ldots),
\end{eqnarray}
and
\begin{eqnarray}\label{48}
\frac{2(n+1)^2}{(n-1)(n^2+2n+10)}<\sin\frac{\pi}{n}<\frac{6(n-1)^2}{(n+1)(n^2-2n+2)}.%\hspace{1.5cm}(n=3,4,\ldots),
\end{eqnarray}
%which are  of the above quantities.
But since
\begin{eqnarray*}
\frac{6(n+1)}{(n+4)(n-2)}-\frac{2(n-1)}{n(n-2)}= \frac{1}{n} \left( \frac{4n^2+8}{n^2+2n-8} \right),
\end{eqnarray*}
\begin{eqnarray*}
\frac{n(n-2)}{n^2-2n+2}-\frac{(n-2)(n+4)}{n^2+2n+10}= \frac{1}{n^2} \left( \frac{16 n^4 -40n^3+16n^2}{n^4+8n^2-16n+20} \right)
\end{eqnarray*}
and
\begin{eqnarray*}
\frac{6(n-1)^2}{(n+1)(n^2-2n+2)}&-&\frac{2(n+1)^2}{(n-1)(n^2+2n+10)}\\
&=& \frac{1}{n} \left( \frac{4n^6-8n^5+44n^4-152n^3+160n^2-64n}{n^6+7n^4-16n^3+12n^2+16n-20} \right),
\end{eqnarray*}
we may write

\begin{eqnarray*}
\tan\frac{\pi}{n}=\frac{2(n-1)}{n(n-2)} +O  \left( \frac{1}{n} \right),
\end{eqnarray*}
\begin{eqnarray*}
\cos\frac{\pi}{n}=\frac{(n-2)(n+4)}{n^2+2n+10} + O \left( \frac{1}{n^2} \right)
\end{eqnarray*}
and
\begin{eqnarray*}
\sin\frac{\pi}{n}=\frac{2(n+1)^2}{(n-1)(n^2+2n+10)} + O  \left( \frac{1}{n} \right).
\end{eqnarray*}
\end{rem}
At the end of this paper we express the following conjecture.\\\\
\textbf{Conjecture.}
It seems that (\ref{46})-(\ref{48}) to be true for all reals $x > 2$ instead of integers $n\geq 3$. The graphs of functions obtained by replacing $n$ by $x$ in (\ref{46})-(\ref{48}) drawn in %corresponding
Figure 1 (a)-(c) respectively strengthen our conjecture. \\
\begin{figure}[h!]
\centering
\subfloat[Related to  (\ref{46})]{\label{fig:tiger}\includegraphics[width=0.30\textwidth]{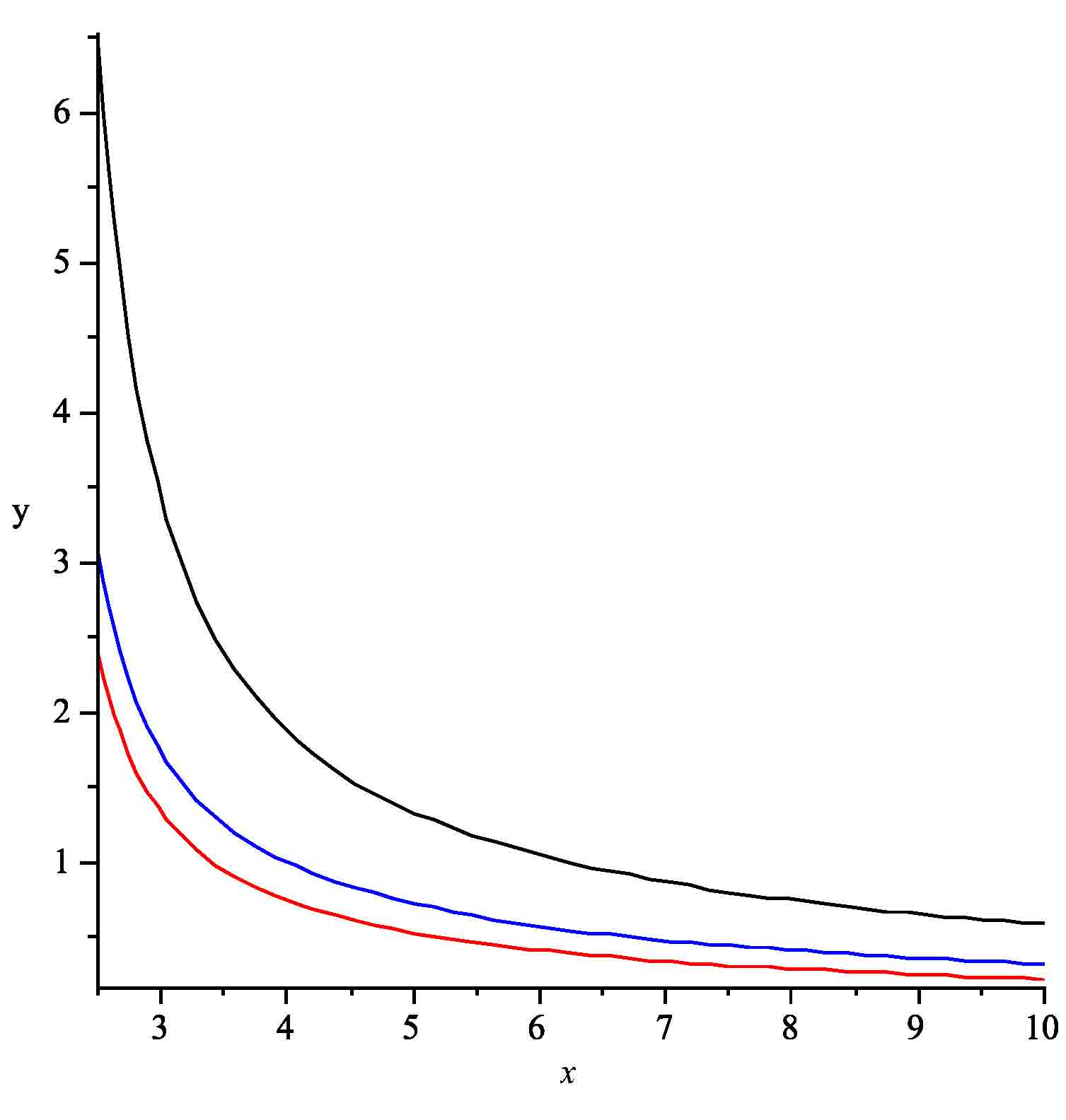}}
  \subfloat[Related to  (\ref{47})]{\label{fig:tiger}\includegraphics[width=0.30\textwidth]{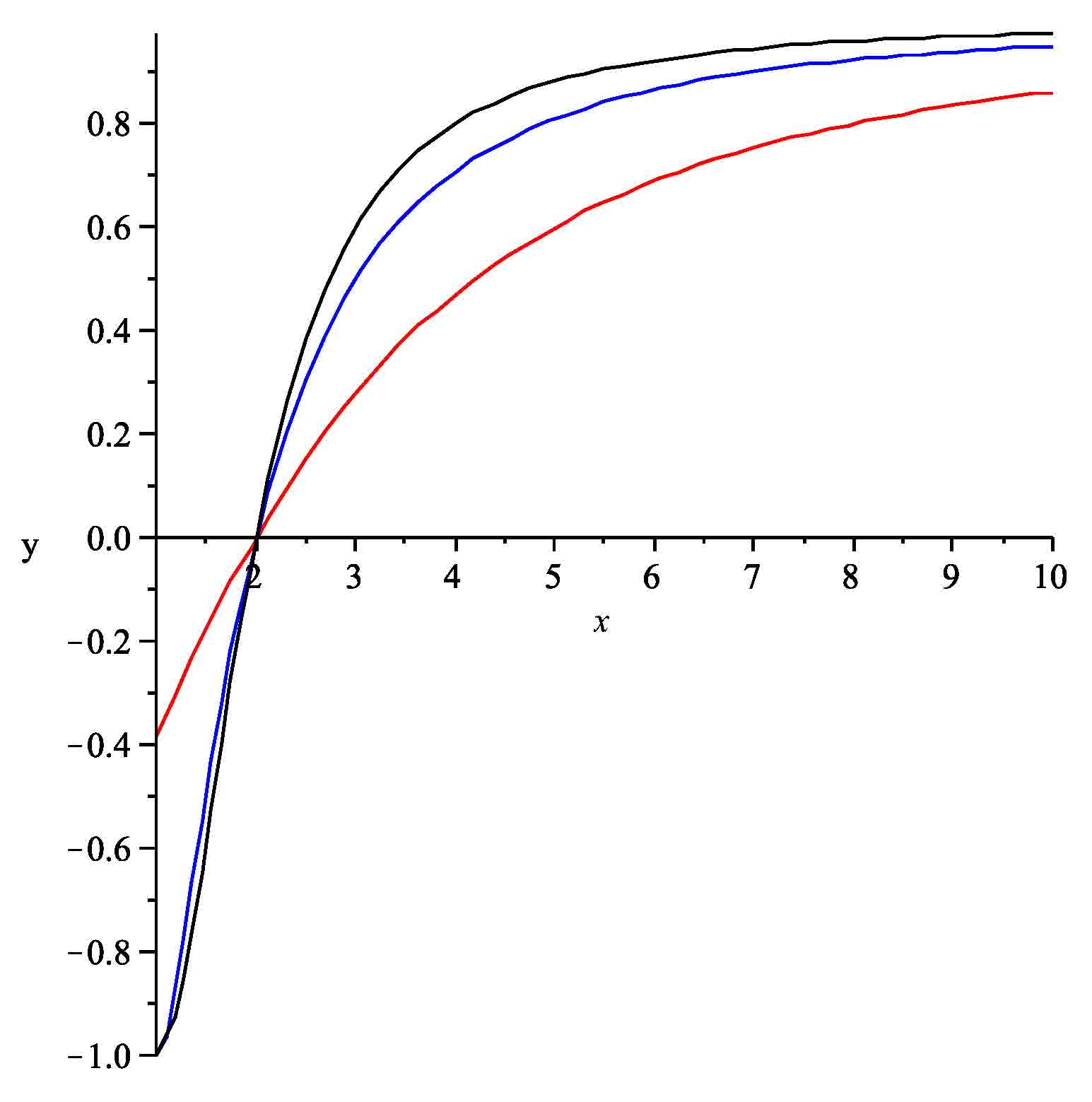}}
  \subfloat[Related to  (\ref{48})]
  {\label{fig:mouse}\includegraphics[width=0.30\textwidth]{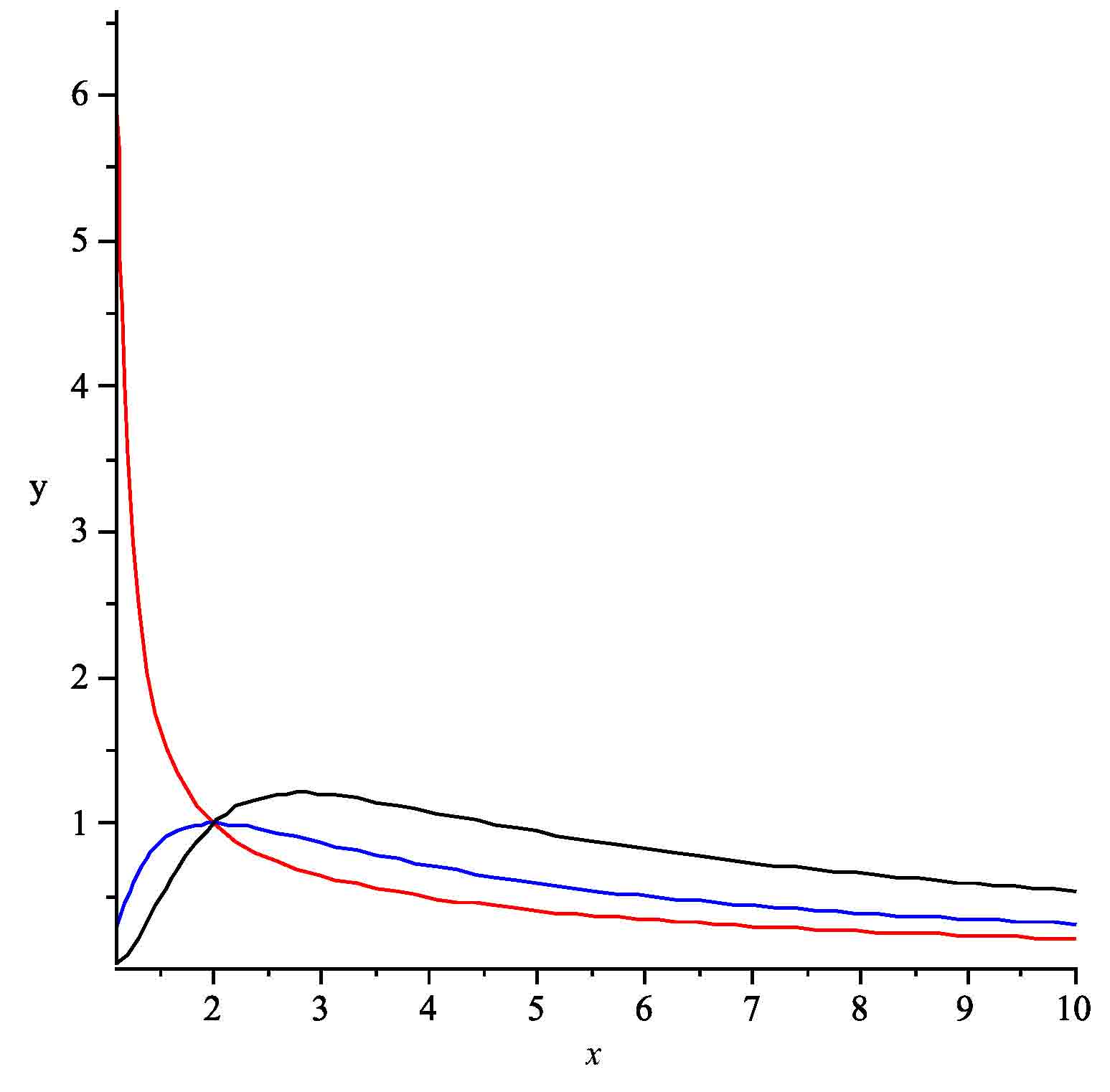}}
\caption{ }
\end{figure}

%\begin{enumerate} \item H. ALZER, {\it An inequality for arithmetic and harmonic means}, Aequations Math. {\bf 46} (1993), 257-263. \item H. ALZER, {\it On an additive analogue of Ky Fan's inequality}, Indag. Mathem., N.S. {\bf 8} (1997), 1-6. \end{enumerate}

\end{document}